\documentclass{article}

\usepackage{amsfonts}
\usepackage{amssymb}
\usepackage{amsthm}
\usepackage[fleqn]{amsmath}
\usepackage[mathcal]{euscript}
\usepackage{mathrsfs}

\newtheorem{coro}{Corollary}

\newtheorem{prop}{Proposition}

\newtheorem{rema}{Remark}
\newtheorem{theo}{Theorem}

\newcommand{\eps}{\epsilon}

\newcommand{\Nset}{\mathbb{N}}
\newcommand{\Rset}{\mathbb{R}}

\newcommand{\dd}{\,{\mathrm{d}}}

\begin{document}


\title{Constructing graphs over $\Rset^n$ with\\ small prescribed mean-curvature}


\author{HOLLY CARLEY \\ 		
\textit{City University of New York}\\ 	
 \\ \sc{and} \\ \\
	MICHAEL	K.-H. KIESSLING \\      
\textit{Rutgers University}\\		
}

\date{Printed: \today}
\maketitle

\begin{abstract}
	In this paper nonlinear Hodge theory and Banach algebra estimates are employed
to construct a convergent series expansion which solves the prescribed mean 
curvature equation $\pm\nabla\cdot(\nabla u/\sqrt{1\pm|\nabla u|^2}) = nH$ for 
$n$-dimensional hypersurfaces in $\Rset^{n+1}$ ($+$ sign) and $\Rset^{1,n}$ ($-$ sign) 
which are graphs $\{(x,u(x)):x\in\Rset^{n}\}$ of a smooth function $u:\Rset^n\to\Rset$,
and whose mean curvature function $H$ 
is $\alpha$-H\"older continuous and integrable, with small norm.
	The radius of convergence is estimated explicitly from below.
	Our approach is inspired by, and applied to, the Maxwell--Born--Infeld theory of electromagnetism 
in $\Rset^{1,3}$, for which our method yields the first systematic way of explicitly computing the electrostatic 
potential $\phi\propto u$ for regular charge densities $\rho\propto H$ and small Born parameter, with explicit
error estimates at any order of truncation of the series.
	In particular, our results level the ground for a controlled computation of Born--Infeld effects 
on the Hydrogen spectrum.
\end{abstract}

\bigskip
\bigskip
\bigskip
\centerline{Revised version: Jan. 24, 2015; Original: Aug. 30, 2010}
\vfill\vfill
\copyright 2015 The authors. Reproduction of this preprint, for non-commercial purposes, is permitted.
\newpage


\section{Introduction}

 In several areas of mathematical physics one encounters the problem of  $n$-dimensional submanifolds 
with prescribed mean-curvature function $H$ in an $n+1$ dimensional Riemannian or Lorentzian manifold.
 For example, the formulation of the Cauchy problem for Einstein's equations in general relativity requires a foliation of the spacetime,
and a popular choice are spacelike so-called constant mean-curvature (CMC) foliations 
\cite{Akutagawa,Bartnik,Gerhardt,Treibergs,MeeksPerez,Fuente};
 in this case the mean curvature function is particularly simple, but entire CMC submanifolds do not always exist and one possible way out 
is to allow the mean curvature to vary in a suitable manner \cite{Gerhardt,MeeksPerezRos}. 
 Another, less obvious example comes from electromagnetism; more precisely, as will be explained in more detail in the next section,
the graph of the electrical potential for the electrostatic solutions 
of the nonlinear Maxwell--Born--Infeld equations in $\Rset^{1,3}$ with prescribed charge density 
\cite{BornA,BornC,Pryce,KieJSPa,CarKie,KieCMP,Gibbons,YangYisong} can be re-interpreted as 
an (almost) entire spacelike hypersurface in Minkowski space \cite{Kly,KlyMik,Bartnik,Gibbons,KieCMP},
with the charge density function playing the role of the mean curvature function.
 Of course, submanifolds with prescribed mean-curvature are also of interest in their own right to differential geometers.

 Although the (non-)existence and regularity theory for solutions of the prescribed mean-curvature equation has evolved quite far
(e.g. \cite{GilbargTrudinger}), the nonlinearity of the equation stands in the way of writing down an explicit solution formula 
except in special circumstances.
 Thus, the problem of constructing an $n$-dimensional hypersurface with explicitly prescribed mean-curvature $H$ in an $n+1$ 
dimensional Euclidean or Minkowskian space has been solved in great generality only when $n=1$ (by simple ODE techniques, 
obviously), although powerful complex variable techniques can be put to work when $n=2$ (especially when 
$H\equiv 0$); see \cite{Courant,Gulliver,Pryce,Kenmotsu,Kobayashi,AkutagawaNishikawa,Struwe,NitscheBOOK,HildebrandBOOKs,
Osserman,OssermanBOOK,Ferraro}.
 However, explicit constructive results when $n\geq 3$ are typically only concerned with radially symmetric solutions (see, e.g., \cite{BJM} 
for a more recent contribution).

 In this paper we are concerned with the explicit construction of solutions to the so-called non-parametric prescribed mean-curvature 
problem where the hypersurfaces are graphs of some scalar function $u$ over the entire $n$-dimensional base manifold when no
symmetry assumption on the prescribed mean-curvature function or its solution is made.
 Yet to find non-special hypersurfaces without symmetry assumptions we have to settle for the neighborhood of
some special solution.
 It is well-known that \emph{in principle} one can systematically find approximate solutions by using 
variational methods, or by running the iterations of the fixed point map of the Schauder theory,
or by directly implementing the proof of the inverse function theorem.
 Curiously enough, though, our literature search for \emph{explicit formulas} of the approximate solutions at arbitrary order, 
with explicit estimates of the domain of convergence of such an approximation scheme, revealed nothing.
 We believe that ours is the first paper in which such a convergent scheme, of explicit quadratures, is presented 
for quite general small and regular mean curvature function; we emphasize once again that no symmetry assumptions are made.

 Specifically, we will explicitly construct a certain neighborhood of the minimal / maximal hypersurfaces ($H\equiv 0$)
which consists of hypersurfaces with small and integrable mean curvature functions $H$ which are graphs over $\Rset^n$, $n \geq 1$.
 In dimension $n=1$, our ``small-curvature'' expansion method actually terminates already after the first term and therefore
yields the solution to the problem by the usual explicit quadrature without smallness assumption.
 Thus our method is primarily of interest in dimensions $n\geq 2$, though for $n=2$ it is perhaps only of interest 
as an alternative to the complex variable techniques.

 Our method is motivated by a concrete problem from physics ($n=3$), more precisely the problem of solving the first-order 
PDE system of Maxwell--Born--Infeld electrostatics. 
 We already mentioned above that the graph of the electrostatic potential is a prescribed mean-curvature hypersurface
in Minkowski space, satisfying the pertinent scalar second-order PDE problem in $n=3$. 
 Therefore in the next section we use this electrostatic analog to motivate and explain our approach, using the 
physicists' familiar three-dimensional vector notation for the fields.
 Then, in section 3, after recalling the scalar prescribed mean-curvature problem for graphs over $\Rset^n$ in all
dimensions $n\geq 1$, in both Euclidean and Minkowskian geometry (our approach treats both on an equal footing),
we use nonlinear Hodge theory to give two dual reformulations of the scalar second-order PDE problem as a first-order 
system for two one-forms which are the $n$-dimensional analogs of the electromagnetic vector equations in $n=3$.
 In section 4 we state our main results in full generality: we present the solutions to 
the first-order Hodge systems, and by corollary to the scalar PDE,
as a convergent series whose terms are given by explicit quadratures.
The parameter of expansion can be seen as a measure of the smallness of $H$ in a suitable H\"older norm.
 In section 5 we prove the convergence in a suitable Banach algebra for each dimension $n\in\Nset$, with $n>1$, 
including an explicit lower estimate of the radius of convergence, obtained with the help of complex analysis; of course 
we also comment on the $n=1$-dimensional case.
 In section 6 we present an algebraically simplified formulation of the three-dimensional problems which
raises some interesting PDE questions.
 We conclude in section 7 by indicating some generalizations of our approach to related quasi-linear PDE problems in divergence form.

\section{The electrostatic analog} 
 
	In Maxwell's classical electromagnetic theory, Coulomb's law states that
an electrostatic charge density $\rho$ in $\Rset^3$ is the source for the so-called 
electric displacement density field $D$,
\begin{equation}
        {\nabla}\cdot D
=
         4 \pi \rho
\phantom{blablabla}\mathrm{Coulomb's\ law}
\,,
\label{eq:COULOMBlaw}
\end{equation}
while Faraday's law says that the so-called electric field strength $E$, when stationary, is 
curl-free:
\begin{equation}
\qquad\qquad \nabla\times{E} = {0} 
\phantom{blablablabla} \mathrm{Faraday's\ law\ (stationary)}
\,.\label{eq:FARADAYlawstatic}
\end{equation}
  The two fields $E$ and $D$ need to be linked by an ``aether law,'' with the help
of which one can eliminate either $D$ or $E$ and obtain a closed, generally nonlinear set of 
first-order PDE for the remaining field.

 Note that ``aether'' is used here merely as a shorthand for ``electromagnetic vacuum;''
historically  Maxwell of course thought of ``\emph{the} classical aether.'' 

 Of particular relevance for the prescribed mean-curvature problem is the aether law
proposed in the 1930s by Max Born \cite{BornA}, which reads
\begin{equation}
D 
= 
\frac{E}{ \sqrt{ 1 - \beta^4 |E|^2 }} 
\phantom{blablabla}\mathrm{Born's\ law}
\phantom{blablablablab}
\label{eq:BORNlawDofE}
\end{equation}
which can be inverted to yield
\begin{equation}
\qquad\qquad E 
= 
\frac{D}{ \sqrt{ 1 + \beta^4 |D|^2 }} 
\phantom{blablabla}\mathrm{Born's\ law\ (reverse)}.
\phantom{blablablablab}
\label{eq:BORNlawEofD}
\end{equation}
 Here, $\beta\in (0,\infty)$ is a hypothetical new constant of nature
(in the dimensionless notation of \cite{KieJSPa}).
 In the limit $\beta\to 0$ Born's law goes over into Maxwell's law of 
the ``pure aether'', $D=E$.

 \begin{rema}: 
{Born's aether law (\ref{eq:BORNlawDofE}) is only the electrostatic
special case of his electromagnetic aether law (replace $|E|^2$ by $|E|^2-|B|^2$,
and add a counterpart for the magnetic induction field $B$ and field strength $H$),
see \cite{BornA,BornC}, and also App. VI in \cite{BornD}.
 Born's electromagnetic aether law was subsequently generalized by Born and Infeld 
\cite{BornInfeldB} to a much more interesting law; the stationary case was recently
treated in \cite{KieJMP}.
 Also higher-dimensional generalizations of Maxwell--Born--Infeld theory with
geometric significance exist; see, for instance, \cite{Gibbons,GibbonsETal,YangYisong}.}
 \end{rema}

  The classically fundamental problem of electrostatic fields generated by point charges, corresponding to
maximal hypersurfaces with conical singularities \cite{Ecker}, is treated in
\cite{BornA,BornC,Pryce,Kly,KlyMik,KieJSPa,CarKie,KieCMP}.
 Unfortunately, only the case with a single point charge \cite{BornC,Ecker}, 
or with a regular lattice of infinitely many of them \cite{Hoppe,Gibbons}, has been solved explicitly. 
 Even the extremely important case of two point charges still awaits solution.

 In this section we make progress by presenting a convergent series solution for electrostatic finite-energy
fields $E$ in all of $\Rset^3$ which are generated by $\alpha$-H\"older continuous charge density functions
$\rho$ which are  sufficiently small in norm. 
 This includes in particular the case of several regularized point charges, i.e. spherically symmetric and compactly
supported $C^\alpha$ functions which we will use in a follow-up work on the old 
--- but still open(!) --- mathematical-physical problem of computing Born--Infeld effects on the quantum-mechanical 
spectrum of atoms, especially Hydrogen.
 This problem is important because the high precision with which the empirical atomic spectra are known will impose
strict bounds on the size of the hypothetical ``new constant of nature'' which enters the Born--Infeld relations ---
thus it will provide a crucial physical viability test of the Born--Infeld proposal for this constant, and possibly on 
the physical viability of their model at that.
 Unfortunately, conflicting results have been obtained so far
with different plausible approximations to the electrostatic pair energy, see \cite{jrlfwg,gsjrwg,gmr,CarKie,Franklin}.
 The solution obtained in this section supplies a series of better and better approximations with controlled error
for smeared out point charges,  which is an important stepping stone along the way to eventually settle the 
problem of Born--Infeld effects on the quantum-mechanical spectra of atoms with point charges.
 Our solution may also pave the ground for computing classical effects in plasmas, cf. \cite{BurtonETal}.

\subsection{Hierarchical vector series solution for small $\beta^2\rho$}

 To state our solution we first stipulate some notation.
 As usual, $C^{k,\alpha}$ denotes the $k$-times continuosly differentiable functions whose $k$-th
derivatives have H\"older exponent $\alpha$; we also write $C^\alpha$ for $C^{0,\alpha}$.
 As for the notation $C_0$, we follow the conventions of \cite{Folland}. 
 Thus, $C_{0}$ denotes the  continuous functions which decay to zero at infinity,
\emph{not} continuous functions with compact support as in \cite{GilbargTrudinger}; the
latter are denoted by $C_c$ in \cite{Folland}, and in \cite{LiebLoss}.
 In the same vein, bounded continuous functions will be denoted by $C_b$ --- here we chose 
not to follow \cite{Folland}, where $BC$ is used instead.
 Altogether, $C^{k,\alpha}_b$, respectively $C^{k,\alpha}_0$, will denote the $k$-times continously 
differentiable functions which together with their derivatives up to order $k$ are bounded, respectively
decay to zero at infinity, and whose $k$-th derivatives are in $C^\alpha$.
 We equip $C^{k,\alpha}_b$ with the usual norm 
\begin{equation}
	\|f\|_{k,\alpha}^{} 
=\label{eq:CkaNORM}
	\max_{|\ell|\leq k}\sup_{x}|D^\ell f(x)|+
	\max_{|\ell|= k}\sup_{x\neq y}\frac{|D^\ell f(x)-D^\ell f(y)|}{|x-y|^\alpha} 
\end{equation}
where $\ell$ is a multi-index.
 Clearly, $C^{k,\alpha}_0\subset C^{k,\alpha}_b$ is a closed subspace, i.e. 
if $f_n^{}\to f$ in $\|\,\cdot\,\|_{k,\alpha}^{}$ norm, and $f_n^{}\in C^{k,\alpha}_0$, then $f\in C^{k,\alpha}_0$, too.

 For $\beta^2\rho\in C^\alpha_{0}\cap L^{1}$ small in the sense stated below we now present
a convergent $\beta$-power series expansion of $D$ solving the system of Maxwell--Born--Infeld equations
(\ref{eq:COULOMBlaw}), (\ref{eq:FARADAYlawstatic}), (\ref{eq:BORNlawEofD}) with vanishing condition at spatial infinity.
 By virtue of (\ref{eq:BORNlawEofD}) this translates into power series solutions for $E$ and 
its electric potential, $\phi$, defined as $E = -\nabla \phi$ with asymptotic vanishing conditions.

\begin{prop}
\label{coro:v}
 Let $\beta^2\rho\in C^\alpha_{0}\cap L^1$ be small in the sense that
\begin{equation}
\label{vGfieldONE}
v^{(1)}(x) := - \nabla \int \frac{\beta^2\rho(y)}{|x-y|}\,\dd^3 y
\end{equation}
has small $C^{1,\alpha}_0$ norm, viz.
\begin{equation}
\label{eq:HsmallnessINvNORM}
\|v^{(1)}\|_{1,\alpha}^{} < \big(2^{2/3}-1\big)^{3/2}.
\end{equation}
 Then the system of PDE 
(\ref{eq:COULOMBlaw}), (\ref{eq:FARADAYlawstatic}), (\ref{eq:BORNlawEofD}) with vanishing condition at spatial infinity
has a solution $D=\beta^{-2} v\in C^{1,\alpha}_0$, with $v$ given by the absolutely convergent series 
\begin{equation}
\label{eq:vSOLUTIONseriesNOeps}
v = \sum_{k=0}^\infty v^{(2k+1)},
\end{equation}
with $v^{(1)}$ given in (\ref{vGfieldONE}) and $v^{(2k+1)}$ for $k\in\Nset$ recursively given by 
\begin{equation}
\label{eq:vsolODDn}
 v^{(2k+1)} =
{\mathbf P} V^{(2k+1)},\qquad k\in\Nset,
\end{equation}
where $V^{(2k+1)}$ is a polynomial in the $v^{(\ell)}$ with odd $\ell<2k+1$, viz.
\begin{equation}
\label{eq:Vdef}
V^{(2k+1)} = 
- \sum_{h=1}^k v^{(2(k-h)+1)}\sum_{j=1}^h M_j^{+}\!
\sum_{|\ell|_{2j} = h-j}\;\prod_{i=1}^jv^{(2\ell_{2i-1}+1)}{\cdot}v^{(2\ell_{2i}+1)},
\end{equation}
with $|\ell|_{K} = \sum\limits_{i=1}^{K} \ell_i$,
the $\ell_i$ taking any non-negative integer values, 
\begin{equation}
\label{eq:Mplus}
M_j^{+} = (- 1)^j\frac{(2j-1)!!}{j!2^j}
\end{equation}
being the $j$-th Maclaurin coefficient of $1/\sqrt{1 +z}$ (with $M_0^+:=1$), and where
 $\mathbf{P}\!: C^{1,\alpha}_{0}\!\to C^{1,\alpha}_{0}$ projects onto the solenoidal 
subspace of $C^{1,\alpha}_{0}$-valued vector fields;
in particular, $\nabla\cdot V^{(2k+1)}\in C^{0,\alpha}_{0}\cap L^{1}$, and so
\begin{equation}
\label{eq:vsolODDnEXPL}
{\mathbf P}V^{(2k+1)}(x) =
V^{(2k+1)}(x) 
+ \nabla \int \frac{{\textstyle{\frac{1}{4\pi}}}\nabla\cdot V^{(2k+1)}(y)}{|x-y|}\,\dd^3 y.
\end{equation}
\end{prop}

 \begin{rema}
 Note that the $v^{(2k+1)}$, $k=0,1,2,...$, depend on $\beta$ in a simple manner, namely
each such $v^{(2k+1)} = \beta^{2(2k+1)}D^{(k)}$, where $D^{(k)}$ is $\beta$-independent.
 Thus, conceptually our solution is nothing but a physicists' formal perturbative series in
powers of $\beta$, treated as a formal small parameter, for which we prove convergence.
 Incidentally, the index $(2k+1)$ rather than $k$ at the $v$s will be explained using perturbation
theory, in section 5.
 \end{rema}

 Having the solution $D=\beta^{-2} v\in C^{1,\alpha}_0$ of the PDE system
(\ref{eq:COULOMBlaw}), (\ref{eq:FARADAYlawstatic}), (\ref{eq:BORNlawEofD}) with vanishing condition at spatial infinity,
we now easily obtain $E$ and $\phi$.

\begin{coro}
\label{coro:wu}
Let $D=\beta^{-2} v\in C^{1,\alpha}_0$ with $v$ given by (\ref{eq:vSOLUTIONseriesNOeps}) solve the PDE system
(\ref{eq:COULOMBlaw}), (\ref{eq:FARADAYlawstatic}), (\ref{eq:BORNlawEofD}) with vanishing condition at spatial infinity.\!
 Then $E=\beta^{-2} w\in C^{1,\alpha}_0$ with $w$ given by 
\begin{equation}
\label{eq:wFROMv}
w = \frac{v}{\sqrt{1 + |v|^2}},
\end{equation}
solves the PDE system
(\ref{eq:COULOMBlaw}), (\ref{eq:FARADAYlawstatic}), (\ref{eq:BORNlawDofE}) with vanishing condition at spatial infinity.

 Moreover, $w$ has its own convergent series expansion, obtained by 
inserting the series solution (\ref{eq:vSOLUTIONseriesNOeps}) for $v$ into the RHS of (\ref{eq:wFROMv})
and expanding, thus
\begin{equation}
\label{eq:powerSERIESw}
w(x) = \sum_{k=0}^\infty w^{(2k+1)}(x),
\end{equation}
with 
\begin{equation}
\label{eq:wPOWERone}
w^{(1)}(x) = v^{(1)}(x),
\end{equation}
while for $k\in\Nset$ each $w^{(2k+1)}$ is given by
\begin{equation}
\label{eq:wPOWERodd}
w^{(2k+1)}(x) 
=
\nabla \int \frac{{\textstyle{\frac{1}{4\pi}}} \nabla\cdot V^{(2k+1)}(y)}{|x-y|}\,\dd^3 y,
\end{equation}
with $V^{(2k+1)}(x)$ given by (\ref{eq:Vdef}) in terms of the $v^{(\ell)}$ with $\ell< 2k+1$.

 Finally, having the solution $E=\beta^{-2}w$ in form of a sum of gradient fields $w^{(2k+1)}$, $k=0,1,2,...$, which
converges in $C^{1,\alpha}_{0}$, we can invoke $E=-\nabla\phi$ and directly read off the electric potential
$\phi=\beta^{-2}u \in C^{2,\alpha}_{0}$, with $u$ given as the series 
\begin{equation}
\label{eq:powerSERIESuASintOFrho}
u(x) =  \int \frac{\rho(y)}{|x-y|}\,\dd^3 y - {\textstyle{\frac{1}{4\pi}}}\sum_{k\in\Nset} 
 \int \frac{\nabla\cdot V^{(2k+1)}(y)}{|x-y|}\,\dd^3 y,
\end{equation}
where $V^{(2k+1)}(x)$ is given by (\ref{eq:Vdef}) in terms of the $v^{(\ell)}$ with $\ell< 2k+1$.
\end{coro}

\begin{rema}
 Note that each term in the expansion of $E = \beta^{-2}w$ is a gradient field, as it should by (\ref{eq:FARADAYlawstatic}).
\end{rema}

\begin{rema}
 The proofs of Proposition 1 and Corollary 1 are given in section 5.
\end{rema}

\begin{rema}
 The expansion technique given here generalizes to electromagnetostatic solutions of the Maxwell--Born--Infeld equations
with regular sources, see \cite{KieJMP}; however, the generalization comes at the price of roughly a factor $1/2$
smaller radius of convergence estimate.
\end{rema}

\subsection{The Helmholtz decomposition}

 As shown by Helmholtz, every vector field in $\Rset^3$ can be decomposed into 
a sum of a divergence-free and a curl-free field. 
 The decomposition is not unique because a harmonic field can be transferred from one 
summand to the other.
 This non-uniqueness can be removed by imposing suitable 
boundary or asymptotic conditions on the individual summands; in our case asymptotic
vanishing conditions.
 The Helmholtz decomposition sheds an interesting light on the $v$ expansion.
 For this reason we include the following elementary observations.

 Writing $v=v_g+v_c$, where $v_g$ is a gradient field (hence, curl-free) 
and $v_c$ is a curl of some vector field (hence, divergence-free), and registering that
${\nabla}\cdot v = {\nabla}\cdot v_g$ and ${\nabla}\times v = {\nabla}\times v_c$, we get
\begin{equation}
\label{vgDIVeqn}
{\nabla}\cdot v_g = 4\pi\beta^2\rho
\end{equation}
(together with $\nabla\times v_g =0$), and 
\begin{equation}
\label{vcCURLeqn}
\nabla\times \frac{(v_g+v_c)}{\sqrt{1\mp |(v_g+v_c)|^2}} =0
\end{equation}
(together with $\nabla\cdot v_c =0$).
 We thereby have obtained a closed set of first-order vector equations for the 
vector field $v_g$ separately from $v_c$, plus a conditionally closed set of 
first-order vector equations for the vector field $v_c$, conditioned on $v_g$ being
given.
 Both sets of PDE are supplemented by the asymptotic conditions that $v_g$ and $v_c$ 
vanish at spatial infinity; with this asymptotic condition the Helmholtz decomposition 
becomes unique.

 As announced, the Helmholtz decomposition sheds an interesting light on the $v$ expansion.
 Since the curl-free component of the Helmholtz decomposition is determined completely
and independently of the divergence-free Helmholtz component $v_c$ of $v$ --- namely thusly:
For $H\in C^\alpha_0\cap L^{1}$ the linear equation (\ref{vgDIVeqn}) has a unique 
gradient field solution $\in C^{1,\alpha}_0\!$ given by
\begin{equation}
\label{vgSOLfield}
v_g(x):= - \nabla \int \frac{\beta^2\rho(y)}{|x-y|}\,\dd^3 y,
\end{equation}
with integration over $\Rset^3$ --- , we can reinterpret our power series approach as a method to solve the $v_c$-equation
(\ref{vcCURLeqn}) given small $v_g$. 
 Indeed, note that $v_g \equiv v^{(1)}$ given earlier in (\ref{vGfieldONE}), so we can just
separate the $v_g$ part off from our power series solution for $v$.
 The remaining power series yields the solution $v_c$ to (\ref{vcCURLeqn}) given small $v_g$. 
 The solution $v_c$ in turn has the meaning of a Born ``correction'' of the Coulomb field $v^{(1)}$  caused by Born's aether law.

 Put in terms of $D$ we summarize: \emph{Our convergent series expansion (\ref{eq:vSOLUTIONseriesNOeps}) yields a solution $D$
given by the $\beta$-power series 
\begin{equation}
D= \sum_{k=0}^\infty\beta^{4k} D^{(k)},
\end{equation}
with $D^{(k)}$ independent of $\beta$, in which $D^{(0)}=\beta^{-2}v^{(1)}$ is the Coulomb field, 
while the $D^{(k)}= \beta^{-2(2k+1)}v^{(2k+1)}$ with $k\in\Nset$ are the coefficient fields of 
its $O(\beta^{4k})$-corrections caused by Born's law.}

\subsection{Prescribed mean-curvature graphs over $\Rset^3$}

 Finally, we also write $w=w_g+w_c$, where $w_g$ is a gradient field (hence, curl-free) and
$w_c$ is a curl of some vector field (hence, divergence-free).
 Clearly, ${\nabla}\cdot w= {\nabla}\cdot w_g$ and ${\nabla}\times w = {\nabla}\times w_c$, and therefore the equation
\begin{equation}
\label{eq:wCURLlaw}
        {\nabla}\times w= 0
\end{equation}
together with the vanishing asymptotics for $w_c$ implies right away that $w_c\equiv 0$, hence $w=w_g=-\nabla u$.
 Now inverting (\ref{eq:wFROMv}) to yield $v=\beta^{-2}D$ in terms of $w=-\nabla u$, and inserting the 
expression in (\ref{eq:COULOMBlaw}) yields a closed second-order scalar equation for $u$, viz.
\begin{equation}
\label{eq:updeMBI}
-\nabla\cdot\frac{\nabla u}{\sqrt{1 -|\nabla u|^2}} = 4\pi\beta^2\rho.
\end{equation}
 This is the well-known equation of a presribed mean-curvature hypersurface in Minkowski spacetime $\Rset^{1,3}$, with
height function $u$ and mean curvature function $H = \frac{ 4\pi}{3}\beta^2\rho$.
 Thus our series expansion for the electrostatic Maxwell--Born(--Infeld) $\phi$ also provides a solution $u$ to the
three-dimensional Minkowskian prescribed mean-curvature equation.

 Replacing $\beta^4\to-\beta^4$ in (\ref{eq:BORNlawDofE}) yields a ``pseudo Born's law,''
in which case the electrostatic Maxwell--Born(--Infeld) equations are
equivalent to the presribed mean-curvature problem in Euclidean space $\Rset^{4}$.
 This equation can therefore be treated with the same expansion strategy applied to an equivalent first order
system of electrostatic Maxwell--pseudo-Born--Infeld equations.

 In fact, our method generalizes to arbitrary dimensions. 
 In the ensuing sections we formulate these problems, state our results, then present their proofs.
 The proofs of our Maxwell--Born--Infeld results are included as special cases.

\section{Formulation of the non-parametric prescribed\\ mean-curvature problems in $\Rset^{n+1}$ and  $\Rset^{1,n}$}
 
 In this section we consider all $n\in\Nset$ equally. 

\subsection{The scalar problem}
	We recall that the mean curvature function $H$ of a hypersurface in $\Rset^{n+1}$
($n=1,2,3,...$), which is the graph $\{(x,u(x)):x\in\Rset^n\}$ of a 
real-valued $C^2$ function $u$ over $\Rset^n$, can be computed from $u$ in a classical manner as
\begin{equation}
H = {\textstyle{\frac1n}}\nabla\cdot\frac{\nabla u}{\sqrt{1+\nabla u\,\cdot\,\nabla u}},
\label{eq:Hpde+}
\end{equation}
(see, e.g., \cite{SerrinA}, and the appendix to section 14 in \cite{GilbargTrudinger}).
	By replacing the Euclidean inner product operation ``$\,\cdot\,$'' by ``$-\,\cdot\,$'' 
one obtains the analogous result for space-like hypersurfaces in Minkowski space $\Rset^{1,n}$, viz.
\begin{equation}
H = -{\textstyle{\frac1n}}\nabla\cdot\frac{\nabla u}{\sqrt{1-\nabla u\,\cdot\,\nabla u}}\,,
\label{eq:Hpde-}
\end{equation}
(see, e.g., \cite{Akutagawa,Bartnik,Kenmotsu,Kobayashi}).
 The inverse problem is to find a function $u$, satisfying 
suitable asymptotic conditions, whose graph over $\Rset^n$ describes a
hypersurface in $\Rset^{n+1}$, respectively $\Rset^{1,n}$, for which a prescribed function 
$H:\Rset^n\to\Rset$ is its mean curvature w.r.t. some normal.
 In that case (\ref{eq:Hpde+}) and  (\ref{eq:Hpde-}) 
become second-order quasi-linear elliptic PDEs for the unknown $u$, 
which we lump together as 
\begin{equation}
\label{eq:HpdePM}
\pm\nabla\cdot\frac{\nabla u}{\sqrt{1\pm|\nabla u|^2}} = nH,
\end{equation}
the upper sign for the Euclidean, the lower for the Minkowskian setting.
	The different signs under the square roots are of course important, while the
overall sign is just a matter of convention; it could be absorbed into $H$.

	Equation (\ref{eq:HpdePM}) needs to be supplemented by an asymptotic condition on $u$ compatible with
the choice of $H$. 
	We will only consider integrable mean curvature functions $H\in L^1$, which implies that $u$ is asymptotic
to some harmonic function, and we choose to impose on $u$ asymptotic monopole conditions, i.e.
\begin{equation}
\label{eq:uHIGHdASYMP}
u(x) \asymp \mp {\textstyle{\frac{n}{|S^{n-1}|}}} \int\! H(y)\dd^n y\, {|x|^{2-n}}
\end{equation}
when $n> 2$, and
\begin{equation}
\label{eq:uTWOdASYMP}
u(x) \asymp
\pm {\textstyle{\frac{1}{\pi}}} \int\! H(y)\dd^2 y\, \ln {|x|}
\end{equation}
when $n=2$.
	However, when $n=1$, then
\begin{equation}
\label{eq:uONEdASYMP}
	u(x) 
\asymp 
\pm \frac{{\textstyle{\frac{1}{2}} \int\! H(y)\dd y}}{\sqrt{1\mp\big|{\textstyle{\frac{1}{2}} \int\! H(y)\dd y}\big|^2}}
|x|;
\end{equation}
clearly, the asymptotic condition (\ref{eq:uONEdASYMP}) with the ``$-$'' sign under the square root
requires the smallness condition
\begin{equation}
\label{eq:ONEdHcondition}
\Big|\int\! H(y)\dd y\Big| <2.
\end{equation}
 This is special to $n=1$ and unrelated to the smallness condition on $H$ which we need to impose later on
to ensure the convergence of an infinite series.

 We note that when $n=1$ or $2$, these asymptotic conditions still allow an arbitrary constant be added to 
any solution $u$ --- we will fix this irrelevant freedom when stating our main results.
	Other options are briefly commented on in the concluding section.

\subsection{The equivalent first-order Hodge systems}

 For all dimensions $n\in\Nset$ the second-order prescribed mean-curvature equation
(\ref{eq:HpdePM}) for either sign is equivalent to a nonlinear Hodge system of first order
for some 1-form.
 There are two mutually dual formulations.

\subsubsection{The $\omega$ system}
 Consider a 1-form $\omega$ satisfying
\begin{equation}
\label{eq:omegaISclosed}
 d\omega =0;
\end{equation}
\begin{equation}
\label{eq:deltaLAWomega}
 \pm \delta \frac{\omega}{\sqrt{1\pm|\omega|^2}} =nH \,.
\end{equation}
 Then this first-order system is equivalent to the second-order equation (\ref{eq:HpdePM}) 
by identifying $\omega \equiv \nabla u \cdot \dd{x}$.
	The asymptotic conditions on $\omega$ are inherited from (\ref{eq:uHIGHdASYMP}), respectively
(\ref{eq:uTWOdASYMP}) or (\ref{eq:uONEdASYMP}); thus, when $n> 2$, 
\begin{equation}
\label{eq:omegaASYMP}
\omega(x) \asymp \mp {\textstyle{\frac{n}{|S^{n-1}|}}} \int\! H(y)\dd^n y\, \nabla{|x|^{2-n}}\cdot \dd{x}
\end{equation}
and when $n=2$, 
\begin{equation}
\label{eq:omegaTWOdASYMP}
\omega(x) \asymp
\pm {\textstyle{\frac{1}{\pi}}} \int\! H(y)\dd^2 y\, \nabla\ln {|x|} \cdot \dd{x}
\end{equation}
while for when $n=1$,
\begin{equation}
\label{eq:omegaONEdASYMP}
\omega(x) \asymp
\pm \frac{{\textstyle{\frac12}\int\! H(y)\dd y}}{\sqrt{1\mp\big|{\textstyle{\frac{1}{2}} \int\! H(y)\dd y}\big|^2}}
\mathrm{sign}(x) \dd{x}.
\end{equation}

\subsubsection{The $\tau$ system}
 Dual to the above system is the following. 
 Defining a 1-form
\begin{equation}
\label{eq:tauDEF}
\tau =  \pm \frac{\omega}{\sqrt{1\pm|\omega|^2}} \,,
\end{equation}
which can be inverted to yield
\begin{equation}
\label{eq:tauDEFinv}
\omega = \pm \frac{\tau}{\sqrt{1\mp|\tau|^2}} \,,
\end{equation}
we see that the above first-order equations are equivalent to
\begin{equation}
\label{eq:deltaLAWtau}
 \delta \tau =nH \,,
\end{equation}
\begin{equation}
\label{eq:dLAWtau}
d \frac{\tau}{\sqrt{1\mp|\tau|^2}} = 0\,.
\end{equation}
 Clearly, (\ref{eq:dLAWtau}) implies that there is a scalar $\sigma$ such that
\begin{equation}
\label{omeageISexact}
\frac{\tau}{\sqrt{1\mp |\tau|^2}} =d\sigma.
\end{equation}
 Up to a sign and an additive constant, $\sigma=u$, of course.

	The asymptotic conditions on $\tau$, inherited from (\ref{eq:uHIGHdASYMP}), respectively
(\ref{eq:uTWOdASYMP}) or (\ref{eq:uONEdASYMP}), are slightly simpler now; namely, when $n=1$,
\begin{equation}
\label{eq:tauONEdASYMP}
\tau(x) \asymp  {\textstyle{\frac{1}{2}}} \int\! H(y)\dd y\, \mathrm{sign}(x)\dd{x}
\end{equation}
and when $n = 2$,
\begin{equation}
\label{eq:tauTWOdASYMP}
\tau(x)
 \asymp
{\textstyle{\frac{2}{|S^{1}|}}} \int\! H(y)\dd^2 y\, \nabla\ln {|x|} \cdot \dd{x},
\end{equation}
whereas for $n>2$, 
\begin{equation}
\label{eq:tauASYMP}
\tau(x) \asymp - {\textstyle{\frac{n}{|S^{n-1}|}}} \int\! H(y)\dd^n y\, \nabla{|x|^{2-n}}\cdot \dd{x}.
\end{equation}

\subsubsection{The cubic $\tau$ system}
 We end this subsection with the observation that upon differentiation in
(\ref{eq:dLAWtau}) and multiplication  with $(1\mp|\tau|^2)^{\frac32}$ we find a cubic version of 
(\ref{eq:dLAWtau}),
\begin{equation} 
\label{eq:tauWEDGElawEXPL}
0 =
\big(1\mp {} |\tau|^2\big) d \tau  
\pm{}d\Big[{\textstyle{\frac12}}  |\tau|^2\Big]\wedge \tau.
\end{equation}
 Together with (\ref{eq:deltaLAWtau}), and supplemented with the same asymptotic conditions, this cubically
nonlinear $\tau$ system is equivalent to the $\tau$ system in the previous subsubsection.
 Ironically, this algebraically simplest formulation of the problem is not at all more user-friendly.
 We shall come back to it in section 6.

\section{Statement of the main results}

 Henceforth we will be primarily concerned with $n\geq 2$.
 As mentioned in the introduction, with some minor adjustments our series solution technique
produces solutions also to the $n=1$-dimensional problem, but the series terminates after the
first term and is identical to the conventional solution by straightforward integration.
 For the details, see our last remark below.

 Our notation of function spaces introduced in section 2 for $n=3$ carries over to arbitrary $n$.
 Now let $H\in C^\alpha_{0}\cap L^1$ and, for $n > 2$, define the exact 1-form
\begin{equation}
\label{tauONEfield}
\tau^{(1)}(x)
 := - 
{\textstyle{\frac{n}{|S^{n-1}|}}}\, d\!\int |x-y|^{2-n}H(y)\,\dd^n y;
\end{equation}
when $n=2$ the logarithmic kernel $\ln\frac{|y|}{|x-y|}$ has to be used in (\ref{tauONEfield}).
 Note that $\tau^{(1)}\in C^{1,\alpha}_0$ for $n\geq 2$. 
 We are now ready to state our main results.
\begin{theo}
\label{thm:tau}
 Let $H\in C^\alpha_{0}\cap L^1$ be small in the sense that
\begin{equation}
\label{eq:Hsmallness}
\|\tau^{(1)}\|_{1,\alpha}^{} < \big(2^{2/3}-1\big)^{3/2}.
\end{equation}
 Then for $n\geq 2$ the Hodge system (\ref{eq:deltaLAWtau}), (\ref{eq:dLAWtau}) with asymptotic condition
(\ref{eq:tauTWOdASYMP}) for $n=2$ and (\ref{eq:tauASYMP}) for $n>2$ has an absolutely convergent 
series solution in $C^{1,\alpha}_0$, given by
\begin{equation}
\label{eq:tauSOLUTIONseriesNOeps}
\tau = \sum_{k=0}^\infty \tau^{(2k+1)},
\end{equation}
with $\tau^{(1)}$ given by (\ref{tauONEfield}) for $n> 2$, and by its logarithmic kernel version for $n=2$,
while $\tau^{(2k+1)}$ for $k\in\Nset$ is recursively given by 
\begin{equation}
\label{eq:tausolODDn}
\tau^{(2k+1)} =
{\mathbf P} T^{(2k+1)},\qquad k\in\Nset
\end{equation}
where $T^{(2k+1)}$ is a polynomial in the $\tau^{(\ell)}$ with odd $\ell<2k+1$, viz.
\begin{equation}
\label{eq:Tdef}
\qquad
T^{(2k+1)} = 
- \sum_{h=1}^k \tau^{(2(k-h)+1)}\sum_{j=1}^h M_j^{\mp}\!\!
\sum_{|\ell|_{2j} = h-j}\prod_{i=1}^j\tau^{(2\ell_{2i-1}+1)}{\cdot}\tau^{(2\ell_{2i}+1)},
\end{equation}
with $|\ell|_{K} = \sum\limits_{i=1}^{K} \ell_i$,
the $\ell_i$ take any non-negative integer values, and with
\begin{equation}
\label{eq:Mpm}
M_j^{\mp} = (\pm 1)^j\frac{(2j-1)!!}{j!2^j},
\end{equation}
the $j$-th Maclaurin coefficient of $1/\sqrt{1\mp z}$ (with $M_0^\mp:=1$), and with
$\tau^{(a)}\cdot\tau^{(b)}$ denoting the conventional inner product of two one-forms
$\tau^{(a)}$ and $\tau^{(b)}$;  moreover, $\mathbf{P}\!: C^{1,\alpha}_{0}\!\to C^{1,\alpha}_{0}$ projects 
onto the co-closed subspace of $C^{1,\alpha}_{0}$(-valued one-forms), i.e. explicitly, for $n > 2$:
\begin{equation}
\label{eq:tausolODDnEXPL}
{\mathbf P}T^{(2k+1)}(x) =
T^{(2k+1)}(x) 
+ {\textstyle{\frac{1}{|S^{n-1}|}}}\, d\! \int \frac{\delta T^{(2k+1)}(y)}{|x-y|^{n-2}}\,\dd^n y,
\end{equation}
while for $n=2$ the logarithmic kernel version has to be used.
\end{theo}

 The series solution for the $\tau$-system immediately yields:
\begin{coro}
	Under the conditions expressed in Theorem \ref{thm:tau}, the 
Hodge system  (\ref{eq:omegaISclosed}), (\ref{eq:deltaLAWomega}) with its pertinent asymptotic condition
(\ref{eq:omegaASYMP}), or (\ref{eq:omegaTWOdASYMP}), has an absolutely convergent 
series solution in $C^{1,\alpha}_0$, which for $n>2$ is given by
\begin{equation}
\label{eq:omegaSOLseriesNOeps}
\omega = 
\mp {\textstyle{\frac{1}{|S^{n-1}|}}} \sum_{k\in\Nset} d \int \frac{\delta T^{(2k-1)}(y)}{|x-y|^{n-2}}\,\dd^n y,
\end{equation}
while for $n=2$ the logarithmic kernel $\ln\frac{|y|}{|x-y|}$ has to be used in (\ref{eq:omegaSOLseriesNOeps}).
 Here, for notational convenience, we have extended the list of $T^{(2k+1)}$s
defined  by (\ref{eq:Tdef}) for $k\in\Nset$ to the case $k=0$, by setting
\begin{equation}
\label{eq:ToneIStauONE}
T^{(1)}:= \tau^{(1)}.
\end{equation}
\end{coro}

 Finally, the series solution for the $\omega$-system now yields:
\begin{coro}
\label{coro:uSERIES}
	Under the conditions expressed in Theorem \ref{thm:tau}, when $n>2$ the scalar PDE (\ref{eq:HpdePM}) 
has an absolutely convergent $C^{2,\alpha}_0$ solution given by
\begin{equation}
\label{eq:uSERIESsol}
u(x) = 
\mp
{\textstyle{\frac{1}{|S^{n-1}|}}} \sum_{k\in\Nset}^\infty 
\int \frac{\delta T^{(2k-1)}(y)}{|x-y|^{n-2}}\,\dd^n y;
\end{equation}
when $n=2$ we use the logarithmic kernel $\ln\frac{|y|}{|x-y|}$, in which case the series converges 
absolutely in $C^{2,\alpha}$ on compact subsets.
\end{coro}

\begin{rema}
The choice of logarithmic kernel in Corollary \ref{coro:uSERIES} fixes the additive constant which
the asymptotic condition (\ref{eq:uTWOdASYMP}) couldn't.
\end{rema}

\begin{rema}
 Alternately, one can of course apply (\ref{eq:tauDEFinv}) 
to the $\tau$-series solution (\ref{eq:tauSOLUTIONseriesNOeps}) 
and map $\tau$ to $\omega$ in order to solve the $\omega$-system.
 However, this is only useful in principle, whereas in practice one will work with the approximate solutions given by
the pertinent partial sums of (\ref{eq:omegaSOLseriesNOeps}).
\end{rema}

\begin{rema}
 When  $n=1$, then the kernel $|y|-|x-y|$ has to be used in (\ref{tauONEfield}), 
and we only have $\tau^{(1)}\in C^{1,\alpha}_b\!$.
 In this case the series terminates right away, i.e. $\tau^{(2k+1)}\equiv 0$ for all $k\in\Nset$,  
and $\tau^{(1)}= \tau$ is the solution.
 However, the ``smallness'' condition (\ref{eq:ONEdHcondition}) has to be imposed in the Euclidean 
setting to map $\tau$ to $\omega$, via (\ref{eq:tauDEFinv}), which is then integrated to get $u$.
\end{rema}

\section{Proofs of the main results}

 We first prove Theorem~1. 
 Setting $n=3$ in the proof proves Proposition~1.

\subsection{Proof of Theorem 1}

The proof consists of two parts: we first construct a formal series solution, and then prove its
absolute convergence in the Banach algebra $C^{1,\alpha}_0$, for which we
obtain an explicit lower estimate for the radius of convergence proportional to the ``size of $H$.''

\subsubsection{The small curvature $\tau$-hierarchy and its formal solution} 

 To facilitate the calculations, we temporarily introduce a ``smallness parameter'' $\eps\in \Rset_+$ through the 
replacement of $H(x)$ by $\eps H(x)$ in (\ref{eq:deltaLAWtau}), then
make the Ansatz
\begin{equation}
\label{eq:powerSERIEStau}
\tau = \sum_{p=1}^\infty \eps^p \tau^{(p)},
\end{equation}
with each $\tau^{(p)}$ independent of $\eps$.
 Inserting (\ref{eq:powerSERIEStau}) into the pair of  equations (\ref{eq:deltaLAWtau}), (\ref{eq:dLAWtau}), 
with $H(x)$ replaced by $\eps H(x)$, and into the asymptotic conditions 
(\ref{eq:tauTWOdASYMP}) for $n=2$ and (\ref{eq:tauASYMP}) for $n>2$, then 
sorting according to powers of $\eps$, we find a hierarchy of linear equations with pertinent asymptotic
conditions.
 At each order in $\eps$ the linear system in question can be solved explicitly in a standard way. 
 At the end of the procedure, we set $\eps=1$ in (\ref{eq:powerSERIEStau}) and
obtain a formal solution of the $\tau$ system  (\ref{eq:deltaLAWtau}), (\ref{eq:dLAWtau}). 

 In particular, $\tau^{(1)}$ satisfies
\begin{equation}
\label{eq:tauDELTAlawONE}
        \delta\tau^{(1)} =nH \,,
\end{equation}
together with
\begin{equation}
\label{eq:tauDlawONE}
d\tau^{(1)} =0.
\end{equation}
 Consistency with our asymptotic conditions for $u$ requires that 
for $H_0\in C^\alpha_{0}\cap L^{1}$ we select the unique solution
$\tau^{(1)}\in C^{1,\alpha}_0$ of the pair of linear equations (\ref{eq:tauDELTAlawONE}), (\ref{eq:tauDlawONE}), 
given by the exact 1-form (\ref{tauONEfield}) for $n> 2$; when $n=2$ our logarithmic kernel will be used.

 Having the exact 1-form (\ref{tauONEfield}), we next find that $\tau^{(2)}$ satisfies
\begin{equation}
\label{eq:tauDELTAlawTWO}
        \delta\tau^{(2)} = 0 
\end{equation}
together with
\begin{equation}
\label{eq:tauDlawTWO}
d\tau^{(2)} =0,
\end{equation}
and $\tau^{(2)} \to 0$ at infinity.
 Clearly, the pair of linear equations (\ref{eq:tauDELTAlawTWO}), (\ref{eq:tauDlawTWO}) 
has a unique solution $\tau^{(2)}$ which tends $\to 0$ at infinity, given 
by the trivial solution $\tau^{(2)}\equiv 0$.

 Carrying on we now find recursively that each $\tau^{(p)}$ for even $p$, i.e. 
$p =2k$ for $k\in\Nset$, satisfies
\begin{eqnarray}
\label{eq:tauDELTAlawEVENord}
        \delta \tau^{(2k)} &=& 0 
\\
\label{eq:tauDlawEVENord}
d\tau^{(2k)} &=&0,
\end{eqnarray}
with $\tau^{(2k)} \to 0$ at infinity,
so that $\tau^{(2k)}\equiv 0$ for general $k\in\Nset$, while for odd $p=2k+1$
with $k\in\Nset$, we find the pair of linear first-order PDE 
\begin{eqnarray}
\label{eq:tauDELTAlawODDk}
\delta\tau^{(2k+1)}&=&0
\\
\label{eq:tauDlawODDk}
d\tau^{(2k+1)}& =&
d T^{(2k+1)},
\end{eqnarray}
where $T^{(2k+1)}$ is the polynomial in the $\tau^{(\ell)}$ with odd $\ell<2k+1$, 
given in (\ref{eq:Tdef}), and ensuing explanations.
 It is easily checked that each $T^{(2k+1)}$ is in $C^{1,\alpha}_{0}$; namely, this
follows inductively from the facts (a) that $\tau^{(1)}\in C^{1,\alpha}_{0}$ when $n\geq 2$, 
and (b) that $C^{1,\alpha}_{0}$ is a Banach algebra.
 Note also that $\delta T^{(2k+1)}\in C^{\alpha}_{0}\cap L^{1}$ for $n\geq 2$, so that (\ref{eq:tausolODDnEXPL})
is well-defined.
 The pair of linear first-order PDE (\ref{eq:tauDELTAlawODDk}), (\ref{eq:tauDlawODDk}) therefore has
a unique solution $\tau^{(2k+1)}\in C^{1,\alpha}_{0}$ given by (\ref{eq:tausolODDn}), with
the projector  $\mathbf{P}\!: C^{1,\alpha}_{0}\!\to C^{1,\alpha}_{0}$ given in (\ref{eq:tausolODDnEXPL})
for $n> 2$; when $n=2$ our logarithmic kernel will  be used.

 Thus we have \emph{formally} solved  our $\tau$-problem (\ref{eq:deltaLAWtau}), (\ref{eq:dLAWtau}), 
with $H$ replaced by $\eps H$, in terms of the formal series solution (\ref{eq:powerSERIEStau}). 
 We can now set $\eps=1$ in (\ref{eq:powerSERIEStau}) and obtain, for $n\geq 2$, a formal series solution of 
(\ref{eq:deltaLAWtau}), (\ref{eq:dLAWtau}) in $C^{1,\alpha}_{0}$, in the sense that 
each partial sum of $\tau= \sum_{k=0}^\infty \tau^{(2k+1)}$ is in $C^{1,\alpha}_{0}$.

 It remains to prove that the so-obtained formal series solution to the $\tau$ system, for $n\geq 2$, 
converges absolutely in $C^{1,\alpha}_{0}$ whenever (\ref{eq:Hsmallness}) holds.

\subsubsection{Absolute convergence of the formal $\tau$ series solution}

 To show that for all $n\geq 2$ the formal series $\tau= \sum_{k=0}^\infty \tau^{(2k+1)}$ converges 
absolutely in $C^{1,\alpha}_{0}$  when the condition (\ref{eq:Hsmallness}) holds, 
it remains to show that the right hand side of the norm estimate
$\|\tau\|_{1,\alpha}\leq \sum_{k=0}^\infty \|\tau^{(2k+1)}\|_{1,\alpha}$ converges.
 To simplify notation we henceforth drop the subscript ``$_{1,\alpha}$'' from the norm symbols.

 We now estimate all norms $\|\tau^{(2k+1)}\|$ for $k\in\Nset$
in terms of the $2k+1$-th power of $\|\tau^{(1)}\|$.
 Now, $\|\tau^{(2k+1)}\| = \|\mathbf{P}T^{(2k+1)}\|$ for $k\in\Nset$, and since 
$\mathbf{P}:C^{1,\alpha}_{0}\to C^{1,\alpha}_{0}$, is a projector, 
we have the estimate: $\|\tau^{(2k+1)}\| \leq \|T^{(2k+1)}\|$.
 Substituting the RHS of (\ref{eq:Tdef}) for $T^{(2k+1)}$, repeating the standard 
inequality $\|\sum_i \tau_i\|\leq \sum_i \|\tau_i\|$, then using the inequality
$\|\tau_i\tau_j\|\leq \|\tau_i\|\|\tau_j\|$ valid in Banach algebras (here
$C^{1,\alpha}_{0}$),  and applying repeatedly the identity
$\|\tau^{(2a+1)}\| = \|\mathbf{P}T^{(2a+1)}\|$ followed by the projector estimate
$\|\mathbf{P}T^{(2a+1)}\|\leq \|T^{(2a+1)}\|$ for the various pertinent values of $a\geq 1$ 
(no estimate is necessary when $a=0$), and using $|M_j^\mp| = M_j^-:=M_j$, 
for $k\geq 1$ we obtain
\begin{eqnarray}
\nonumber
\|T^{(2k+1)}\| &\leq &
\sum_{h=1}^k \|\tau^{(2(k-h)+1)}\|\sum_{j=1}^h |M_j^\mp|\!\!
\sum_{|\ell|_{2j} = h-j}\prod_{i=1}^{2j}\|\tau^{(2\ell_i+1)}\|
\\
\label{eq:TnormESTIM}
&\leq&
\sum_{h=1}^k \|T^{(2(k-h)+1)}\|\sum_{j=1}^h M_j
\sum_{|\ell|_{2j} = h-j}\prod_{i=1}^{2j}\|T^{(2\ell_i+1)}\|.
\end{eqnarray}

 Next we show that for all $k\in\Nset$ there exists some $R_{2k+1}$ such that
\begin{equation}
\label{eq:TnormESTIMiterated}
\qquad
\|T^{(2k+1)}\|\leq  R_{2k+1} \|\tau^{(1)}\|^{2k+1}.
\end{equation}
 Setting $k=1$ and recalling the definition (\ref{eq:ToneIStauONE}) we obtain the estimate
\begin{equation}
\label{eq:TnormESTIMone}
\qquad
\|T^{(3)}\| \leq  {\textstyle{\frac12}} \|\tau^{(1)}\|^{3}.
\end{equation}
 Now suppose that for all $k=1,...,k_*$ there exists some $R_{2k+1}$ such that (\ref{eq:TnormESTIMiterated}) holds.
 Then the estimate (\ref{eq:TnormESTIM}) guarantees that (\ref{eq:TnormESTIMiterated})
is true also for $k=k_*+1$, and since $k_*\geq 1$ is arbitrary in this induction step while
(\ref{eq:TnormESTIMone}) says that the estimate is true for $k_*=1$, it follows that
(\ref{eq:TnormESTIMiterated}) is true for all $k\in\Nset$.

 The inductive proof that (\ref{eq:TnormESTIMiterated}) holds for all $k\in\Nset$
also yields that $R_{2k+1}$ is recursively defined for $k\in\Nset$ by
\begin{equation}
\label{eq:Rrecursion}
\quad R_{2k+1} = 
\sum_{h=1}^k R_{2(k-h)+1}\sum_{j=1}^h M_j
\sum_{|\ell|_{2j} = h-j}\prod_{i=1}^{2j}R_{2\ell_i+1},
\end{equation}
with $R_1:=1$. 

 Recall that we want to show that $\sum_{k=0}^\infty  \|\tau^{(2k+1)}\|<\infty$ 
for sufficiently small $\|\tau^{(1)}\|$.
 Since we have shown that $\|\tau^{(2k+1)}\| \leq \|T^{(2k+1)}\| \leq R_{2k+1} \|\tau^{(1)}\|^{2k+1}$, 
it suffices to show that 
$\sum_{k=0}^\infty R_{2k+1} \|\tau^{(1)}\|^{2k+1}< \infty$ for sufficiently small $\|\tau^{(1)}\|$.

 Setting now $\|\tau^{(1)}\|=:\xi$, we note that the formal power series
$G(\xi) := \sum_{k=0}^\infty R_{2k+1}\xi^{2k+1}$ is nothing but the formal generating function 
of the $R_{2k+1}$, in the usual sense that, formally, $R_{2k+1} = G^{({2k+1})}(0)/({2k+1})!$.
 So our task is to show that the generating function is analytic about $\xi=0$ with
radius of convergence given by the RHS of (\ref{eq:Hsmallness}).
 
 With the help of the recursion relation (\ref{eq:Rrecursion}) we readily find that $G(\xi)$ 
is the positive inverse function of $g\mapsto \xi$ given by
\begin{equation}
\label{eq:GimplicitDEF}
\xi = 2g -  \frac{g}{\sqrt{1-g^2}}
\end{equation}
defined for positive $\xi$ near $\xi=0$, with $G(0)=0$.
 Since the function $g \mapsto \xi$ given by (\ref{eq:GimplicitDEF}) is analytic about 
$g=0$ (with radius of convergence $=1$) and has unit derivative at $g=0$, there now 
exists an open neighborhood of $\xi=0$ on which there is defined a unique inverse function 
$\xi\mapsto g=G(\xi)$ which vanishes at $\xi=0$, has unit derivative at $\xi=0$, and 
satisfies (\ref{eq:GimplicitDEF}).
 Thus, in particular, the Maclaurin expansion of $G(\xi)$ converges to $G(\xi)$ and it
generates the recursion coefficients $R_{2k+1}$.

 We now determine the finite radius of convergence $\xi_*$ of the power series for $G(\xi)$ 
about $\xi=0$.
 Setting $G(\xi) = \sin \Psi(\xi)$ we see that (\ref{eq:GimplicitDEF}) defines
the function $\psi\mapsto \xi$ given by
\begin{equation}
\label{eq:sinGimplicitDEF}
\xi = 2\sin\psi  - \tan\psi,
\end{equation}
with $\xi=0$ when $\psi=0$. 
 Since $\psi\mapsto \sin\psi$ is an entire function, which vanishes for $\psi=0$ and
has unit derivative there, the radius of convergence of the Maclaurin series of 
$\xi\mapsto G(\xi)= \sin\Psi(\xi)$ coincides with the radius of convergence of 
the Maclaurin series of $\xi\mapsto \Psi(\xi)$. 
 This radius of convergence in turn is determined by those $\xi$ value(s) closest to 
$\xi=0$ at which the derivative of $\psi\mapsto \xi = 2\sin\psi -\tan\psi$ vanishes 
(possibly asymptotically should $\xi\to\xi_\infty$ when $|\psi|\to\infty$ suitably).
 But this $\psi$ derivative is $2\cos\psi - 1/\cos^2\psi$, and it vanishes iff $2\cos^3\psi=1$,
which gives $2^{1/3}\cos\psi\in\{1,e^{i2\pi/3},e^{i4\pi/3}\}$.
 A calculation now gives the radius of convergence of $G(\xi)$ about $\xi=0$ as
\begin{equation}
\label{eq:radOFconv}
 \xi_* = \big(2^{2/3}-1\big)^{3/2}.
\end{equation}

 Thus we have shown that our formal series solution converges absolutely
 in $C^{1,\alpha}_{0}$ if $\|\tau^{(1)}\| < \big(2^{2/3}-1\big)^{3/2}$.
 This completes our convergence proof.
\qed

\begin{rema}
 Because of our use of the projector and Banach algebra estimates,  we cannot conclude that
$\|\tau^{(1)}\| < \big(2^{2/3}-1\big)^{3/2}$ is a necessary criterion for absolute convergence
of our formal power series solution.
 Indeed, for any radially symmetric $H\in C^\alpha_0$ the formal power series reduces to
its first term, all other terms being identically zero, so convergence is a trivial issue
and holds for any size of $H$, then.
\end{rema}

\subsection{Proofs of Corollaries 2 and 3}

 To prove Corollary 2, use (\ref{eq:tauDEFinv}) but insert (\ref{eq:powerSERIEStau}) as well as
the analogous Ansatz
\begin{equation}
\label{eq:powerSERIESomega}
\omega = \sum_{p=1}^\infty \eps^p \omega^{(p)},
\end{equation}
then sort by powers of $\eps$ and easily solve for $\omega^{(p)}$ explicitly; reinsert those expressions
into (\ref{eq:powerSERIESomega}) and now set $\eps=1$. 
 The result is (\ref{eq:omegaSOLseriesNOeps}).
 Its absolute convergence follows verbatim as in the proof of Theorem 1.

 To prove Corollary 3, integrate $\omega = \nabla u\cdot dx$, with $\omega$ from given by the
absolutely convergent series (\ref{eq:omegaSOLseriesNOeps}); the constant of integration is fixed by 
the asymptotic conditions.
\qed 

\newpage
\section{The cubic reformulation of the $\tau$-problem}

 We now come back to our earlier remark that the nonlinearity of the $\tau$-problem is effectively cubic.
 However, for technical reasons, in this section we restrict attention to $n=3$.
 This allows us to return to the three-dimensional vector formulation for which many well-known
identities of vector analysis are at our disposal.
 Thus the $\tau$-problem then becomes a $v$-problem.

 The ``cubic version'' of our $v$-problem reveals some interesting 
a priori differential identities which are satisfied by any solution but which 
are obscured by the original formulation of the $v$-problem.

 More importantly, the algebraic simplifications dramatically reduce the 
combinatorial complexity of the successive approximations in the small-$H$ hierarchy.
 Indeed, except for the first few low-order terms, the series solution to our original 
first-order $v$ vector problem and its spin-off, the solutions to the original
first-order $w$ vector problem and the original second-order scalar $u$-problem, 
soon involves terms which look more and more unwieldy. 
 The cubic version of the $v$-problem offers relief, by having to evaluate fewer integrals.

 Interestingly enough, though, the solution theory of this apparently simpler cubical formulation
of the small-$H$ hierarchy seems more complicated than the ``square-root formulation'' and 
raises some challenging questions.
 More to the point, it will be clear from the equivalence of the original and the cubic $v$-problems that the 
cubic version has a solution $v$ whenever $v$ solves the original problem, and vice-versa. 
 Moreover, upon reintroducing the $\eps$ parameter, it is clear from the analyticity in $\eps$ 
that a small-$H$ expansion will produce the same series solutions as before. 
 However, if we ignore for a moment that we already know that our formal series derived from the
original version of the $v$-problem converges absolutely to a classical solution, from
which the solvability of the cubic hierarchy of linear equations follows as a corollary,
then the  solvability of the cubic hierarchy of linear equations is not at all obvious but needs to be verified!
 We will prove the consistency ab initio, but did not succeed in proving its convergence without recourse to the
square-root formulation.
 Perhaps some reader will feel inspired to settle this problem!

\subsection{\!\!The cubically nonlinear reformulation of the $v$-problem}

 Carrying out the curl operation in 
\begin{equation}
\label{eq:vCURLlaw}
\nabla\times \frac{v}{\sqrt{1\mp |v|^2}} =0
\end{equation}
and multiplying through with $(1\mp|v|^2)^{3/2}$ we find that, away from singularities (which occur, for instance, 
when $1-|v|^2 \to 0$), any solution $v$ of (\ref{eq:vCURLlaw}) satisfies 
\begin{equation} 
\label{eq:vCURLlawEXPL}
0 =
\big(1\mp {} |v|^2\big) \nabla\times v   \pm{}\Big[{\textstyle{\frac12}}\nabla |v|^2\Big]\times v.
\end{equation}
 This equation already exhibits a cubic nonlinearity, yet
it can be further manipulated into a more concise alternate format.

 Namely, using first the identity 
${\textstyle{\frac12}}\nabla |v|^2 = (v\cdot\nabla)v+v\times(\nabla\times v)$
well-known from vector analysis, we find that (\ref{eq:vCURLlawEXPL}) is equivalent to
\begin{equation} 
\label{eq:vCURLlawEXPLrewr}
0 =
\big(1\mp {} |v|^2\big) \nabla\times v
  \pm{}[(v\cdot\nabla)v+v\times(\nabla\times v)]\times v.
\end{equation}
 We next employ the identity
$v\times(\nabla\times v)\times v =|v|^2\nabla\times v -v(v\cdot \nabla\times v)$,
well-known from vector algebra, which for any solution of (\ref{eq:vCURLlawEXPL}) simplifies 
to $v\times(\nabla\times v)\times v = |v|^2\nabla\times v$ because any solution of 
(\ref{eq:vCURLlawEXPL}) satisfies (\ref{curlVdotV}).
 Inserting $v\times(\nabla\times v)\times v = |v|^2\nabla\times v$ into
(\ref{eq:vCURLlawEXPLrewr}) yields a cancellation and (\ref{eq:vCURLlawEXPLrewr}) becomes
\begin{equation}
\label{veqnROT}
\nabla\times v
=
\pm{}v\times (v\cdot\nabla)v.
\end{equation}
 This is perhaps the most concise cubic form of the curl equation for $v$.
 When paired with the divergence equation 
\begin{equation}
\label{eq:vDIVlaw}
        {\nabla}\cdot v =3H \,
\end{equation}
we arrive at the following conclusion.

\begin{prop}
\label{prop:cubicv}
 The curl equation (\ref{veqnROT}) paired with the divergence equation (\ref{eq:vDIVlaw}) 
forms a closed system of first-order vector PDE for $v$.
 To make them well posed the vector equations need to be supplemented by asymptotic conditions
for $v$ at spatial infinity, which as before we take to be vanishing in agreement with the 
monopole asymptotics of $u$ at spatial infinity.
 By construction, this cubic set of equations for $v$, with $v$ vanishing at infinity, is equivalent to the scalar 
equation (\ref{eq:HpdePM}), with $u$ vanishing at infinity, for the appropriate choice of sign. 
\end{prop}

\subsection{The Helmholtz decomposition for the cubic version}
 
 For the sake of completeness, we remark 
that the Helmholtz decomposition $v=v_g+v_c$ leads to analogous conclusions for the cubic 
version of the $v$-problem.
 We find the previously obtained linear equation (\ref{vgDIVeqn}) for $v_g$
(together with $\nabla\times v_g =0$), while the nonlinear equation for $v_c$ (given $v_g$)
becomes
\begin{equation}
\label{vReqn}
	\nabla\times v_c=\pm{}(v_g+v_c)\times (v_g+v_c)\cdot\nabla(v_g+v_c)
\end{equation}
(together with $\nabla\cdot v_c =0$).
 This closed set of first-order vector equations for the 
vector fields $v_g$ and $v_c$ is supplemented by the asymptotic conditions
that $v_g$ and $v_c$ vanish at spatial infinity.

 As announced, given the gradient field (\ref{vgSOLfield}), the remaining equation
(\ref{vReqn}) now becomes a closed vector equation for the solenoidal field $v_c$.
 It remains to solve equation (\ref{vReqn}) with $v_g$  given by (\ref{vgSOLfield}).

\subsection{Spin-off: a-priori differential identities for solutions}
\noindent
 The cubic version of the $v$-problem reveals two a-priori differential identities 
which are satisfied by any solution of the $v$-problem.
\begin{prop}
\label{prop:vIDs}
 For any solution of the $v$-problem,
\begin{equation}
\label{curlVdotV}
v\cdot(\nabla\times v) = 0
\end{equation} 
and
\begin{equation}
\label{nabvsqrdotcurlv}
\nabla |v|^2  \cdot \nabla\times v 
=
0,
\end{equation}
wherever $\nabla\times v$ is defined.
\end{prop}

\noindent
\textit{Proof of Proposition \ref{prop:vIDs}:}

\noindent
 Dotting (\ref{eq:vCURLlawEXPL}) with $v$ yields (\ref{curlVdotV}), and 
dotting (\ref{eq:vCURLlawEXPL}) with $\nabla{|v|^2}$ yields (\ref{nabvsqrdotcurlv}). \qed

 Clearly, neither (\ref{curlVdotV}) nor (\ref{nabvsqrdotcurlv}) 
are generally true for arbitrary vector fields $v$.

\subsection{The small-$H$ expansion: cubic version}

 We now turn to the \emph{formal} small-$H$ solution strategy for the cubic $v$-problem.
 We begin by deriving the cubic analog of the hierarchy of linear PDE.
 To present a slightly different perspective, this time we take the 
Helmholtz point of view.
 Also, we take the liberty and bypass the reintroduction of $\eps$ and, instead of 
``powers of $\eps$,'' simply talk about ``orders of smallness.'' 

 In this vein, for $H\in C^\alpha_{0}\cap L^{1}$ we find $v_g$ given in (\ref{vgSOLfield}).
 Now, if $H$ is suitably small, then $v_g$ is small in $C^{1,\alpha}_{0}$, and 
we write $v_g = v^{(1)}$.
 For the solenoidal part of $v$ we  now make the series Ansatz
$v_c= v^{(p_1)} + v^{(p_2)}+ v^{(p_3)} + \dots$ with
$1< p_1<p_2<p_3<\cdots$, and with each $v^{(p_k)}$ being of $p_k$-th ``order of smallness,'' compared to $v^{(1)}$.
 Inserting this Ansatz into the cubic vector PDE (\ref{vReqn}) and identifying $p_1$ with 
the smallest order on the RHS, $p_2$ with the next-to-smallest order, and so on, we find 
recursively that $p_1=3$, then $p_2=5$, and $p_k=2k+1$ for general $k\in\Nset$.
 Furthermore, for $v = v_g +v_c$ with $v_g = v^{(1)}$ and 
$v_c= v^{(3)} + v^{(5)}+  v^{(7)} + \dots$
to be a solution of the pair of equations (\ref{eq:vDIVlaw}), (\ref{vReqn}), 
each $v^{(2k+1)}$ for $k\in\Nset$ has to satisfy 
\begin{equation}
\label{curlEQorderK}
\nabla\times v^{(2k+1)} =
\pm \sum_{h+i+j =  k-1} v^{(2h+1)}\times (v^{(2i+1)}\cdot\nabla)v^{(2j+1)}
\end{equation}
supplemented by the solenoidality condition 
\begin{equation}
\label{divEQorderK} 
\nabla\cdot v^{(2k+1)}=0.
\end{equation}
 Supposing that $v^{(2\ell+1)}$ is known for all $\ell\leq k-1$, then 
(\ref{curlEQorderK}), (\ref{divEQorderK})
is a pair of linear first-order PDE for $v^{(2k+1)}$, with vanishing
conditions at spatial infinity for $v^{(2k+1)}$.
 Now $v^{(1)}$ is known, and so, by induction, it follows that 
(\ref{curlEQorderK}), (\ref{divEQorderK}) successively determine
$v^{(3)}$, then $v^{(5)}$, and so on --- provided that each equation in this formal
linear hierarchy of equations is solvable (in $C^{1,\alpha}_{0}$)!

\subsection{Solvability of the cubic hierarchy}

\begin{prop}
\label{prop:solvabilityCUBIC}
The infinite hierarchy of equations (\ref{curlEQorderK}), (\ref{divEQorderK}), 
together with $v_g = v^{(1)}$ given in (\ref{vgSOLfield}), is uniquely solvable in 
$C^{1,\alpha}_{0}\cap L^2$ at each order $k\in\Nset$, and this solution is given by
\begin{equation}
\label{solORDERk}
v^{(2k+1)}(x)
=
\pm{\textstyle{\frac{1}{4\pi}}}\sum_{h+i+j\atop =k-1} 
\int (v^{(2h+1)}\times (v^{(2i+1)}\cdot\nabla)v^{(2j+1)})(y)\times\frac{x-y\ }{|x-y|^3}\dd^3y.
\end{equation}
\end{prop}

\noindent
\textit{Proof of Proposition \ref{prop:solvabilityCUBIC}:}

\noindent
 Since
$\sum_{h+i+j=k-1}v^{(2h+1)}\times(v^{(2i+1)}\cdot\nabla)v^{(2j+1)}\in C^{\alpha}_{0}\cap L^{1}$
if $v^{(2\ell+1)}\in C^{1,\alpha}_{0}\cap L^{2}$ for all $\ell\leq k-1$, it
remains to be verified that the collected terms on the RHS of (\ref{curlEQorderK}) 
have vanishing divergence at each order $k$.
 With some effort one can show term by term that (\ref{solORDERk}) 
is identical to (\ref{eq:vsolODDn}), (\ref{eq:vsolODDnEXPL}). 
 For instance, based on the identity $\nabla(a\cdot b) = (a\cdot\nabla) b + (b\cdot\nabla)a
+a\times\nabla\times b + b \times\nabla\times a$ for any two vector fields $a$ and $b$ it
is readily shown that 
\begin{equation}
-v^{(1)}\times (v^{(1)}\cdot\nabla)v^{(1)}
=
\frac{1}{2}\nabla|v^{(1)}|^2\times v^{(1)}
=
\nabla\times({\textstyle{\frac12}}|v^{(1)}|^2v^{(1)}),
\end{equation}
proving the equality of the RHS of (\ref{solORDERk}) 
and that of (\ref{eq:vsolODDn}), (\ref{eq:vsolODDnEXPL}) when $k=1$.
 However, the procedure of proving equality term by term  soon gets very complicated, 
and an inductive argument is needed, instead.

 Proceeding by induction, we will now show that, given $v_g = v^{(1)}$, for each $k\in\Nset$ 
there is a $v^{(2k+1)}$ which vanishes at infinity and solves
\begin{equation}
\label{ih}
\nabla\times v^{(2k+1)} =
\pm\sum_{h+i+j=k-1} v^{(2h+1)}\times (v^{(2i+1)}\cdot\nabla)v^{(2j+1)}.
\end{equation}
 Since the formulas become rather long, we switch to the shorter
notation $v^{[k]}$ for $v^{(2k+1)}$ for all $k=0,1,2,3,...$;
this mildly obscures the order of smallness to which the terms belong, but
shortens the length of the formulas considerably.

 Step one is easily disposed of by checking explicitly that
$v_g\times (v_g\cdot\nabla) v_g$ is divergence-free for any gradient field $v_g$.
 Indeed, since $\nabla\times v_g =0$ for any gradient field $v_g$, 
first of all the well-known identity 
$(v\cdot\nabla)v = {\textstyle{\frac12}}\nabla |v|^2 - v\times(\nabla\times v)$ 
reduces to
$(v_g\cdot\nabla)v_g = {\textstyle{\frac12}}\nabla |v_g|^2$, 
and second,
$v_g\times \Big({\textstyle{\frac12}}\nabla |v_g|^2\Big)
= 
-\nabla\times \Big[{\textstyle{\frac12}}|v_g|^2v_g\Big]$,
so that
$v_g\times (v_g\cdot\nabla) v_g = 
-\nabla\times \Big[{\textstyle{\frac12}}|v_g|^2v_g\Big]$ is a curl,
i.e. divergence-free.
 Hence, given $v^{[0]} = v_g$, for either sign of ``$\pm$'' there is a 
$v^{[1]}$ satisfying $\nabla\times v^{[1]}= \pm v^{[0]}\times (v^{[0]}\cdot\nabla) v^{[0]}$.
 Moreover, since $v_g$ is given by (\ref{vgSOLfield}) and $H\in C^\alpha_{0}\cap L^{1}$, 
we have that $v_g\times (v_g\cdot\nabla) v_g\in C^{1,\alpha}_{0}\cap L^{1}$, and so in 
particular there is a solution 
$v^{[1]}$ of $\nabla\times v^{[1]}= \pm v^{[0]}\times (v^{[0]}\cdot\nabla) v^{[0]}$
which vanishes at infinity.

 As to the induction step, suppose that for all $k\leq m$, and either sign of 
``$\pm$,'' there is a solution $v^{[k]}$ of (\ref{ih}) which vanishes at spatial 
infinity.
 We now show that then also
$\sum_{h+i+j=m} v^{[h]}\times (v^{[i]}\cdot\nabla)v^{[j]}$ is divergence free
and vanishes sufficiently rapidly at infinity.



 We will need two equalities.
 Taking the dot product of both sides 
of (\ref{ih}) with $v^{[l]}$ and summing for $l+k=m'$  with $m'\leq m$ gives 
\begin{eqnarray}
\sum_{k+l=m'}{v^{[l]}}\cdot \nabla\times v^{[k]}
 &=&
\pm \sum_{k+l\atop =m'}\sum_{h+i+j\atop =k-1} 
{v^{[l]}}\cdot \Big[v^{[h]}\times (v^{[i]}\cdot\nabla)v^{[j]}\Big]
\\
&=& 
\pm \sum_{h+i+j+l\atop =m'-1}
{v^{[l]}}\cdot \Big[v^{[h]}\times (v^{[i]}\cdot\nabla)v^{[j]}\Big]\\
&=&
 0\label{eq1},
\end{eqnarray}  
where the last equality follows by the anti-symmetry of the triple product
$a\cdot(b\times c)$;  (\ref{eq1}) is the order $m$ 
expansion analogue of (\ref{curlVdotV}).

 Similarly, dotting both sides of (\ref{ih}) with 
$(v^{[l]}\cdot\nabla)v^{[p]}$ and summing over $k+l+p=m$ gives
\begin{eqnarray}
&&\sum_{k+l+n  =m}\nabla\times v^{[k]}\cdot((v^{[l]}\cdot\nabla)v^{[p]})=
\\
&&\sum_{k+l+n\atop  =m}\sum_{h+i+j\atop =k-1} 
(v^{[h]}\times (v^{[i]}\cdot\nabla)v^{[j]})\cdot ((v^{[l]}\cdot\nabla)v^{[p]})=\\
&&
\sum_{h+i+j+l+n  =m-1}
 (v^{[h]}\times (v^{[i]}\cdot\nabla)v^{[j]})\cdot ((v^{[l]}\cdot\nabla)v^{[p]})=0.
\label{eq2}
\end{eqnarray}


 Now note that
\begin{eqnarray*}
&&
\sum_{h+i+j  =m}v^{[h]}\times (v^{[i]}\cdot\nabla) v^{[j]} =
\\
&&
\sum_{h+i+j  =m}
v^{[h]}\times \Big({\textstyle{\frac12}}\nabla(v^{[i]}\cdot v^{[j]})
-v^{[i]}\times(\nabla\times v^{[j]})\Big) =\\
&&
\sum_{h+i+j  =m}\Big(-\nabla\times
\big({\textstyle{\frac12}} (v^{[i]}\cdot v^{[j]}) v^{[h]}\big)+
{\textstyle{\frac12}}(v^{[i]}\cdot v^{[j]})
\nabla\times v^{[h]}\\
&&
\hskip2truecm
-(v^{[h]}\cdot \nabla\times v^{[j]})v^{[i]}
+(v^{[h]}\cdot v^{[i]}) \nabla\times v^{[j]}\Big)=\\
&&
\sum_{h+i+j  =m}
\Big(-\nabla\times({\textstyle{\frac12}} (v^{[i]}\cdot v^{[j]}) v^{[h]})
+{\textstyle{\frac32}}(v^{[h]}\cdot v^{[i]})\nabla\times v^{[j]}\Big),
\end{eqnarray*}
where we used that 
$\sum_{h+i+j=m} (v^{[h]}\cdot \nabla\times v^{[j]})v^{[i]}=0$, by (\ref{eq1}). 
 And so,
\begin{eqnarray*}
\quad&&\nabla\cdot\!\!\sum_{h+i+j  =m}\!\! v^{[h]}\times (v^{[i]}\cdot\nabla)v^{[j]}
=
{\textstyle{\frac32}}\sum_{h+i+j  =m} \nabla(v^{[h]}\cdot v^{[i]})\cdot\nabla\times v^{[j]}
=\\
&&\sum_{h+i+j  =m}
\Big[\Big(
(v^{[h]}\cdot\nabla)v^{[i]} + (v^{[i]}\cdot\nabla)v^{[h]}\Big)\cdot\nabla \times v^{[j]}\\
&&\hskip1.5truecm 
+\Big((v^{[h]}\times\nabla \times v^{[i]}+
v^{[i]}\times \nabla\times v^{[h]})\Big)\cdot\nabla \times v^{[j]}\Big].
\end{eqnarray*}
 By (\ref{eq2}), 
$\sum_{h+i+j=m} \big((v^{[h]}\cdot\nabla)v^{[i]} +(v^{[i]}\cdot\nabla)v^{[h]}\big)
\cdot\nabla \times v^{[j]}=0$.
 Moreover, by the anti-symmetry of the triple product $a\cdot(b\times c)$ we have
\begin{equation}
\sum_{h+i+j=m} 
(v^{[h]}\times(\nabla \times v^{[i]})+v^{[i]}\times(\nabla\times v^{[h]}))\cdot\nabla \times v^{[j]}
=0. 
\end{equation}
 Thus $\sum_{h+i+j=m} v^{[h]}\times (v^{[i]}\cdot\nabla)v^{[j]}$ is divergence free,
as claimed. 
\qed

\subsection{Convergence of the cubic $v$ series}

 Summing (\ref{solORDERk}) over $k\in\Nset$ and adding $v^{(1)} = v_g$ given by (\ref{vgSOLfield}) 
yields a formal series solution to the cubic $v$-problem, viz.
\begin{eqnarray}
v(x)
=\nonumber
\hskip-.5truecm
&&
 v_g(x) \pm \\
\label{formalSERIESsol}
&&{\textstyle{\frac{1}{4\pi}}}
\sum_{k=1}^\infty \!\!
\sum_{h+i+j=\atop k-1}\!\! \nabla\times\int 
\frac{(v^{(2h+1)}\times (v^{(2i+1)}\cdot\nabla)v^{(2j+1)})(y)}{|x-y|}\dd^3y.
\end{eqnarray}
 We now have to address the convergence of this formal series solution.

 Since the combinatorial structure of the RHS of (\ref{solORDERk}) is considerably simpler than
that of (\ref{eq:vsolODDn}), (\ref{eq:vsolODDnEXPL}), it would seem that a convergence proof
is more readily forthcoming than our previous proof. 
 Curiously, the obviously simpler series expansion for the cubic version of the $v$-problem 
has not at all yielded to our attempts of proving its convergence \emph{directly}, 
i.e. without recourse to the original version of the first-order 
vector problem with its more complicated nonlinearity. 
  So we are finally forced to recall the origin of the cubic version of the $v$-problem to
conclude that (\ref{formalSERIESsol}) converges  absolutely to a
classical solution for small $H$.
\medskip

 So there is also a challenge: it would be good to have a simple convergence
proof directly for the cubic hierarchy!

\section{Related problems in divergence form}

 Our power series technique of solving the prescribed mean-curvature problem for graphs over 
$\Rset^n$ can handle more general quasilinear problems in divergence form. 
 In particular, let $f(s)= 1 + as + bs^2 + \cdots$ be given and analytic about $s=0$, then
the equation 
\begin{equation}
\nabla\cdot \Big(f(|\nabla u|) {\nabla u}\Big) = \rho,
\label{eq:RHOpde}
\end{equation}
with $\rho\in C^\alpha_0\cap L^1$ small enough and asymptotic monopole condition for $u$,
can be solved with the same solution techniques developed here for the
scalar prescribed mean-curvature equation.
 In particular, such type of problems occur in the theory of stationary compressible fluid 
flows, 
see \cite{Bers,SiSiB,Shiffman,Smith}, and other versions of nonlinear electrostatics \cite{Demetri}.

 While we here are interested in entire solutions over $\Rset^n$ motivated by the electromagnetic
Born--Infeld model, and also by some problems in the theory of spacetime structure 
(e.g. \cite{Akutagawa,Bartnik,Gerhardt,Treibergs}), with some mild modifications our approach should 
be adaptable to the prescribed mean-curvature equation in bounded domains with small Dirichlet 
data for $u$  (e.g. \cite{SerrinA,BombieriGiusti,Giusti,Trudinger,Takakuwa,Bourni,COOR}); 
this then also includes the minimal / maximal hypersurface 
problem ($H\equiv 0$) with small Dirichlet data 
(e.g. \cite{Courant,NitscheBOOK,OssermanBOOK,HildebrandBOOKs,Struwe,Obersnel,Bergner}). 
 The parametric prescribed mean-curvature problem
(e.g. \cite{Wente,BrezisCoron,CaldiroliMusina,Struwe,Ye,ChanilloMalchiodi}), 
which also captures embedded hypersurfaces which are not graphs over $\Rset^n$ and
even non-embedded hypersurfaces, has a  different structure, however.

\smallskip


\textbf{Acknowledgement} This work was started while H. Carley was a Hill Assistant Professor at
Rutgers University, supported in parts by the NSF through grant DMS-0406951.
      M. Kiessling was supported by the NSF through grants DMS-0406951 and
DMS-0807705; and in the final stages of this work also by CNRS Universite de Provence
through a poste rouge. 
 Thanks go to A. Shadi Tahvildar-Zadeh and Sagun Chanillo for helpful comments on
the manuscript, and to Y. Elskens for his hospitality and interest in this work.

 \vfill\newpage

        



\frenchspacing
\bibliographystyle{plain}

\end{document}